\documentclass[11pt,notitlepage,oneside]{amsart}

\usepackage{amssymb}
\usepackage{amsthm}
\usepackage{booktabs}

\theoremstyle{plain}

\newtheorem{thm}{Theorem}

\theoremstyle{definition}

\begin{document}

\title{On a class of reducible trinomials}

\begin{abstract}
In this short note we give an expression for some numbers $n$
such that the polynomial $x^{2p}-nx^p+1$ is reducible.
\end{abstract}

\maketitle

\begin{center}
Ralf Stephan\footnote{\texttt{ralf@ark.in-berlin.de}}
\end{center}


{\bf Subject Class:\ }

\bigskip

Trinomial reducibility was generally treated by
Schinzel\cite{schinzel2007selecta} but he only mentions in passing
in his first theorem the case of a degree~2 factor.
Filaseta et al\cite{filaseta2007two} give some easy criteria
for reducibility of
some trinomials. Bremner\cite{bremner1981trinominals} determines
all trinomials $x^n+Ax^m+1$ with irreducible cubic factor. In a more
recent work\cite{bremner2011reducibility}, Bremner and Ulas find
reducibility criteria for several trinomials; in particular they show
in theorem~4.4 the conditions for factorization of $x^6+Ax^3+B$ into
three factors of degree~2.

Not many polynomials where only one coefficient is variable allow for a
general treatment of their irreducibility. Let $P(m,A)=x^{2m}-Ax^m+1$,
then trivially $A=2$ makes $P$ reducible for all $m$. Also, by
substitution, if $P(m,A)$ is reducible then so is $P(km,A),k>1$.
The interesting cases are therefore the polynomials $P(p,A)$, with
odd prime $p$.

In summary we have
\begin{thm}
The polynomial $P(p,A)=x^{2p}-Ax^p+1$ is reducible if
\begin{equation}\label{f1}
A=\sum_{0\le i\le(p-1)/2}(-1)^{\left(i+\frac{p-1}{2}\right)}
\frac{p}{2i+1}\binom{i+\frac{p-1}{2}}{2i}k^{2i+1},
\quad k=1,2,3\dots 
\end{equation}
In particular,
\begin{equation}\label{f2}
x^{2p}-Ax^p+1=(x^2-kx+1)\cdot Q(k,p),
\end{equation}
with $Q$ a palindromic polynomial having coefficients from a subset
of the values of the Lucas-type sequence defined by $a_{i+2}=ka_{i+1}-a_i$,
$a_1=0$, $a_2=1$.
\end{thm}
An example would be $p=5$, $k=3$:
\[
x^{10}-123x^5+1=(x^2-3x+1)(x^8+3x^7+8x^6+21x^5+55x^4+21x^3+8x^2+3x+1).
\]

Integer sequences $a_i$ that satisfy linear recurrences with constant
coefficients have the property that
$c_ia_i+c_{i+1}a_{i+1}+\ldots+c_{i+h}a_{i+h} = 0$
for some $h$, and this can be used to construct two polynomials
which when multiplied vanish at nearly all coefficients. This paper
is concerned with $h=2$.

Let $a_i$ a Lucas-type sequence defined by $a_{i+2}=ka_{i+1}-a_i$, or
$a_i-ka_{i+1}+a_{i+2}=0$. Let $s_i, 1\le i\le2p+1$ a finite
palindromic sequence with $s_i=s_{2p+2-i}$ and for $i\le p+1$, $s_i=a_i$.
The term $s_i-ks_{i+1}+s_{i+2}$ is zero for all $i$ except when
$i=p+1$, and the process of computing the term $s_i-ks_{i+1}+s_{i+2}$
over all $1\le i\le2p-1$ is equivalent to multiplying $\sum_i s_ix^{i-1}$
(the polynomial in $x$ with coefficients $s_i$)
by the polynomial $x^2-kx+1$, the resulting nonzero terms being
$x^{2p}$, $x^p$, and 1.

It remains to derive the coefficient $A(k,p)$ of $x^p$ which equals
$ks_{p+2}-2s_{p+1}$. This is equivalent to
\begin{align*}
A(k,p)\quad=\quad &ka_{p+1}-2a_{p}\\
\quad=\quad &a_{p+2}-a_{p}\\
\quad=\quad &\displaystyle{[z^p]\frac{1-z^2}{1-kz+z^2}}\\
\quad=\quad &\displaystyle{[z^p]\left(\frac{1-z^2}
{1+z^2}\cdot\frac{1}{1-\frac{z}{1+z^2}k}\right)}\\
\quad=\quad &\sum_{0<i\le p} d_{i,p} k^{p-i-1},\\
\end{align*}
with $d$ the elements of the Riordan array $\mathcal{R}\left(\frac{1-z^2}
{1+z^2},\frac{z}{1+z^2}\right)$, from which the proposition
follows (see for example Merlini\cite{b-mra}).
As example, $x^{10}-Ax^5+1$ is reducible if $A$ is of form $k^5-5k^3+5k$.

\end{document}